# Learning to Configure Mathematical Programming Solvers by Mathematical Programming [*]

Gabriele Iommazzo[†]    Claudia D'Ambrosio[†]
dambrosio@lix.polytechnique.fr

Antonio Frangioni[‡]
frangio@di.unipi.it

Leo Liberti[†]
liberti@lix.polytechnique.fr

January 8, 2024

**Abstract**

We discuss the issue of finding a good mathematical programming solver configuration for a particular instance of a given problem, and we propose a two-phase approach to solve it. In the first phase we learn the relationships between the instance, the configuration and the performance of the configured solver on the given instance. A specific difficulty of learning a good solver configuration is that parameter settings may not all be independent; this requires enforcing (hard) constraints, something that many widely used supervised learning methods cannot natively achieve. We tackle this issue in the second phase of our approach, where we use the learnt information to construct and solve an optimization problem having an explicit representation of the dependency/consistency constraints on the configuration parameter settings. We discuss computational results for two different instantiations of this approach on a unit commitment problem arising in the short-term planning of hydro valleys. We use logistic regression as the supervised learning methodology and consider CPLEX as the solver of interest.

## 1 Introduction

Mathematical Programming (MP) is a formal language for describing optimization problems; once a problem is modelled by a MP formulation, an off-the-shelf solver can be used to solve it. Off-the-shelf solvers must be general enough to encompass a significant family of problems, and yet fast enough that sufficiently large-scale instances will be solved in reasonable time. By the usual trade-off between generality and efficiency, implementing a good solver is extremely hard. Today's most successful solvers, such as e.g. IBM-ILOG CPLEX [17], meet these specifications by actually embodying a corpus of different solution algorithms, each with their own (often very large) set of algorithmic options [18]. The default values for these options are usually chosen so that the solver will perform reasonably well on a large instance library, but for many problem classes performances can be improved by carefully tuning the solver configuration. Identifying good solver parameters for a given problem instance, however, is a difficult art, which requires a considerable experience in solver usage, and an in-depth hands-on knowledge of the application giving rise to the considered MP formulation.

---

[*]This paper has received funding from the European Union's Horizon 2020 research and innovation programme under the Marie Skłodowska-Curie grant agreement n. 764759 "MINOA".

[†]LIX CNRS, École Polytechnique, Institut Polytechnique de Paris, Palaiseau, France. Gabriele Iommazzo is now at Zuse Institute Berlin, Berlin, Germany (iommazzo@zib.de).

[‡]Dip. di Informatica, Università di Pisa, Pisa, Italy

Automatic configuration of algorithmic parameters is an area of active research, going back to the foundational work in [26]. Many approaches are based on sampling algorithmic parameter values, testing the solver performance, and performing local searches in order to find the parameters that most improve performance [14, 9, 16, 22]. An algorithmic configuration method, derived from [16] and specifically targeted to MP solvers (including CPLEX), is described in [15]. All of these methods learn from a set of instances the best configuration for a similar set of instances; the configuration provided is not "per-instance" but common to all of the instances in the same problem class. This is different to the approach investigated in the present paper, which aims at providing a specific configuration for each given instance. The per-instance approach is necessary whenever the solver performance on instances of a problem class varies much, as in the case of our specific application. A more theoretical approach to choosing provably optimal parameter values based on black-box complexity, and limited to evolutionary algorithms, is given in [8] and references therein. Many artificial intelligence methodologies have been applied to this problem, see e.g. [23].

While the previously cited methodologies — and the one proposed in this paper — try to learn the best parameter values of an algorithm before launching it to solve a new instance, other approches learn on-the-fly, during its execution. For instance, the CPLEX Automatic Tuning Tool [17, Ch. 10], accompanying the corresponding solver, runs it on an instance (or a set thereof) several times, within a user-decided time limit, testing a specific parameter setting at each run, and saves the configuration providing the best algorithmic performance. Another on-the-fly methodology is presented, e.g., in [3], where one or more parameters of a tabu-search heuristic are adjusted, during its execution and in function of its behaviour, until a good configuration is learned.

In this work we present a new two-phase approach to automatic solver configuration. In the first phase, called the *Performance Map Learning Problem* (PMLP), we use a supervised Machine Learning (ML) methodology in order to learn the relationships between the features $f$ of an instance $\iota$, the solver configuration $c$ and the performance quality $p$ of the solver configured by $c$ on the instance $\iota$. Formally, $p$ is defined as a function $p : (f, c) \to \mathbb{R}$ measuring the integrality gap achieved by the solver within a certain time limit. We propose two different variants of the PMLP. In the Performance-as-Output (PaO) one, we learn an approximation $\bar{p}(f, c)$ of the performance function. In the Performance-as-Input (PaI) one, we instead learn an approximation $\bar{\mathcal{C}}$ of the map $\mathcal{C} : (f, r) \to c$, where $c$ is any configuration allowing to obtain a required performance level $r \in \mathbb{R}$ for the instance $\iota$. In the second phase we use use the learnt information (either $\bar{p}$ or $\bar{\mathcal{C}}$) to define the *Configuration Space Search Problem* (CSSP), a constrained optimization problem, which is different for the two variants PaO and PaI. The input of the CSSP is an encoding of the performance map as well as the features of the instance to be solved. Its constraints encode the logical dependency and compatibility conditions over the configuration parameters: this ensures feasibility of the produced configuration. The objective function of the CSSP is a proxy to the performance map: optimizing it over the the constraints yields a good configuration for the given instance.

Our approach is therefore capable of handling configuration spaces having arbitrarily complex logical conditions. This overcomes a weakness in previous learning-based approaches to the algorithm configuration problem, as acknowledged e.g. by [15, 16]. To see how this weakness might adversely impact a solver configuration methodology, consider the following naive approach: learn the map $\mathcal{C}$ using a supervised ML method, then ask the trained method to output $\mathcal{C}(f, 1)$ (1 being the best possible performance) for a new, unseen instance encoded by $f$. Unfortunately, this approach would fail over most off-the-shelf supervised ML methodologies, which are unable to reliably enforce dependency and compatibility constraints on the output configuration. Some attempts have been made to overcome this issue. For instance, the authors of ParamILS [16], that performs local searches in configuration space, declare that their algorithm supports the encoding of dependence/compatibility constraints on feasible parameters configurations. Unfortunately, we were unable to find the precise details of this encoding. Other approaches, used in learning-based optimization try to directly integrate a constrained optimization problem in a neural network, embedding it into the gradient computations of the back-propagation pass [11, 29] or into an individual layer [2]. However, they are not generalizable to any ML algorithm and/or MP. We make one last introductory remark about the parameter search space: obviously, our approach can help configure any subset of solver parameters; in order to reduce the time spent in constructing the training set, a judicious choice would consider a reasonably small subset of parameters that are thought to have a definite impact on the problem at hand.

The rest of this paper is organized as follows. In Sect. 2 we formally introduce the notation and the main ingredients of our approach. In Sect. 3 we discuss both variants of the PMLP, and in Sect. 4 we discuss the corresponding CSSP. Finally, in Sect. 5 we report computational experiments and we draw some conclusions.

## 2 Notation and preliminaries

### 2.1 The training set

Let $C$ be the set of valid solver configurations. We assume for simplicity that $C \subseteq \{0, 1\}^s$, although extension to integer and continuous numerical parameters is clearly possible. Since every subset of the unitary hypercube can be described by means of a polytope [25, Cor. 1], we assume that its representation as a set of linear inequalities in binary variables, say

$$C = \{\, c \in \{\,0,\,1\,\}^s \mid Ac \le d \,\} \tag{1}$$

is known. In practice, deriving $A$ and $d$ from the logical conditions on the parameters can be assumed to be easy.

Let $\mathcal{O}$ be an optimization problem consisting of an infinite number of instances. In order to be able to use a ML approach, we have to encode each instance $\iota \in \mathcal{O}$ by a feature vector $f_\iota \in \mathbb{R}^t$ for some fixed $t \in \mathbb{N}$. This is surely possible at least by restricting ourselves to some subset of $\mathcal{O}' \subseteq \mathcal{O}$ (say, instances with appropriately bounded size). We also assume availability of a finite subset $I \subset \mathcal{O}'$ of instances and let $F = \{\, f_\iota \mid \iota \in I \,\}$ be their feature encodings. We remark that $F$ must be representative of $I$. Since feature extraction is an intensively studied field, we do not dwell on the specifics here. We also remark that ML methodologies are known to perform well on training sets that are not "overly general" [13, Ch. 5.3]: thus, we assume that $I$ is a set of instances belonging to the same problem, or at least to different variants of a single problem.

In practice, in this paper we focus on a unit commitment problem arising in the energy industry, which is solved hundreds of times per day. The instances all have the same size; the constraints do not vary overmuch; the features are the objective function coefficients. Notwithstanding, our approach is general: the "problem structure" is encoded in the set of features (extracted from the instances), which are certainly class-specific, but need not be size-specific (one can e.g. use dimensionality reduction techniques to achieve feature vectors of the same size even for instances of varying size).

We can then, in principle, compute $p(f, c)$ on each feature vector $f \in F$ with every configuration $c \in C$ by calling the solver configured with $c$ on the instance $\iota$, in order to exploit $(F \times C, p(F, C))$ as a training set for learning estimates $\bar{p}$ or $\bar{\mathcal{C}}$ as described above. Hopefully, then, these can be used to generalize our approach to instances outside $I$, with known encoding and that are in some way similar (in size or otherwise) to those in I.

### 2.2 Logistic regression

Logistic Regression (LR) is a supervised ML methodology devised for binary classification of vectors [7].

Let $\mathcal{X} = (\mathcal{X}_1, \ldots, \mathcal{X}_m)$ be a vector of random variables, and let $\mathcal{Y}$ be a Bernoulli distributed random variable depending on $\mathcal{X}$. Following [20], and denoting $\mathsf{P}(\mathcal{X} = x)$ by $\mathsf{P}(x)$ and $\mathsf{P}(\mathcal{Y} = y)$ by $\mathsf{P}(y)$, we have

$$\begin{aligned}
\mathsf{P}(1|x) &= \frac{\mathsf{P}(x|1)\mathsf{P}(1)}{\mathsf{P}(x)} = \frac{\mathsf{P}(x|1)\mathsf{P}(1)}{\mathsf{P}(x|1)\mathsf{P}(1) + \mathsf{P}(x|0)\mathsf{P}(0)} \\
&= \frac{1}{1 + \frac{\mathsf{P}(x|0)\mathsf{P}(0)}{\mathsf{P}(x|1)\mathsf{P}(1)}} = \frac{1}{1 + e^{-z}} = \sigma(z),
\end{aligned} \tag{2}$$

$$\text{where } z = \ln\frac{\mathsf{P}(x|1)}{\mathsf{P}(x|0)} + \ln\frac{\mathsf{P}(1)}{\mathsf{P}(0)}. \tag{3}$$

We now assume that $z$ depends linearly on $x$:

$$\exists w \in \mathbb{R}^m, b \in \mathbb{R} \quad z = wx + b. \tag{4}$$

In some cases $w, b$ can be computed explicitly: for example, if we assume that the conditional probabilities $\mathsf{P}(x|y)$ are multivariate Gaussians with means $\mu_y$ and identical covariance matrices $\Sigma$, and use the above expression for $\mathsf{P}(1|x)$, we obtain

$$\mathsf{P}(1|x) = \frac{1}{1 + e^{-(wx+b)}}$$

where $w = \Sigma^{-1}(\mu_1 - \mu_0)$ and $b = \frac{1}{2}(\mu_0 + \mu_1)^\top \Sigma^{-1}(\mu_0 - \mu_1) + \ln(\mathsf{P}(1)/\mathsf{P}(0))$. In general, however, explicit formulæ cannot always be given, and $w, b$ must be computed from sampled data.

To simplify notation in this section, we let $\tau(x, w, b) \triangleq \frac{1}{1+e^{-(wx+b)}}$, simply denoted $\tau(x)$ when $w, b$ are fixed. Since we only consider two class labels $\{0, 1\}$, we model the probability of $\mathcal{Y} = y$ conditional to $\mathcal{X} = x$ using the function $\tau(x)$ if $y = 1$ and $1 - \tau(x)$ if $y = 0$, i.e.
$$P(y \mid x) = \tau(x)^y (1 - \tau(x))^{1-y}, \quad (5)$$
which evaluates to $\tau(x)$ whenever $y = 1$ and $1 - \tau(x)$ whenever $y = 0$.

We consider a training set $T = (X, Y)$ where $X = (x^i \in \mathbb{R}^m \mid i \leq n)$ and $Y = (y^i \in \{0, 1\} \mid i \leq n)$ consist of independent and identically distributed samples. We then use the Maximum Likelihood Estimation methodology [13] to find the optimal values of the parameters $w$ and $b$. To this end, we define the likelihood function
$$L_T(w, b) = \prod_{i \leq n} P(y^i \mid x^i) = \prod_{i \leq n} \tau(x^i, w, b)^{y^i} (1 - \tau(x^i, w, b))^{1-y^i}. \quad (6)$$
We want to maximize $L_T(w, b)$. Since the logarithmic function is monotone in its argument, maximizing $\ln(L_T(w, b))$ yields the same optima $(w^*, b^*)$ as maximizing $L_T(w, b)$. The training problem of LR is therefore:
$$\max_{w,b} \sum_{i \leq n} \left[ y^i \ln \left( \frac{1}{1 + e^{-(wx^i + b)}} \right) + (1 - y^i) \ln \left( 1 - \frac{1}{1 + e^{-(wx^i + b)}} \right) \right] \quad (7)$$
We recall that the functions
$$\psi_1(z) = \ln \left( \frac{1}{1 + e^{-z}} \right) \quad \text{and} \quad \psi_2(z) = \ln \left( 1 - \frac{1}{1 + e^{-z}} \right)$$
are concave [6, Ex. 3.49(a)]. Since $0 \leq y^i \leq 1$ for each $i \leq n$, Eq. (7) maximizes the sum of convex combinations of concave functions, so it is a convex optimization problem which can be solved efficiently.

Once trained, the LR maps input vectors $x \in \mathbb{R}^m$ to an output scalar $y \in [0, 1]$: in this sense, LR approximates a binary scalar; binary values can of course be retrieved by rounding, if necessary. We denote this by $y = \mathsf{LR}(x, w^*, b^*)$.

### 2.3 The performance function

Since any LR output must be in $[0, 1]$ by definition, the performance data set $p(F, C)$ must be scaled to $[0, 1]$. We first measured, on different instance and solver configuration and within a given time limit, the CPLEX integrality gap, which is defined as
$$\frac{|\text{best integer sol.value} - \text{best relaxation value}|}{1\mathrm{e}{-10} + |\text{best integer sol.value}|}$$
in [17], pg. 263. Unfortunately, CPLEX performance data sometimes include very large values which stand for "infinity" (denoted $\infty$ below), meaning that CPLEX may not find feasible solutions or valid bounds at every run within the allotted CPU time. Instead of scaling in the presence of these outliers, we employed the following algorithm to obtain our performance function:

1. fix a constant $\gamma > 0$;
2. let $\hat{p} = \max(p(F, C) \smallsetminus \{\infty\})$;
3. for each $v \in p(F, C)$ if $v > \hat{p}$ let $v = \hat{p} + \gamma$;
4. let $\rho = (\rho_1, \ldots, \rho_n)$ be a ranking of $p(F, C) = (v_1, \ldots, v_n)$ so that $\rho_1 \geq \rho_2 \geq \cdots \geq \rho_n$ and $\rho_1$ ranks the best performance value (equal values in $p(F, C)$ are assigned the same rank);
5. scale $\rho$ to $[0, 1]$ so that $\rho_1$ is mapped to $1$ and $\rho_n$ to $0$.

The choice of LR for this work is motivated by the fact that: (a) the parameters chosen for automatic configuration are all binary, and LR is a good method for estimating binary values; (b) the performance function has range $[0, 1]$. In general, LR can be replaced by other ML methodologies: this changes the technical details of the two phases, but it does not change the overall approach.

## 3 The PMLP

We now describe in details the two announced variants of the PMLP.

### 3.1 PMLP-PaO

In this variant, the output that we want to produce is an approximation $\bar{p}(f, c)$ of the performance function. Therefore, we interpret the symbols in Sect. 2.2 using the entities defined in Sect. 2.1. We note that the $y$ variables in Eq. (7) are continuous, as noted at the end of Sect. 2.2. We have $X = F \times C \subseteq \mathbb{R}^{t+s}$ and $Y = \rho$; that is, $x = (f, c)$ in (7) encodes the concatenation of features and configurations, and $y = \rho$ is a vector of dimension 1 (i.e., a scalar) encoding the performance (see Fig. 1 (left)). We optimize (7) using some local Nonlinear Programming (NLP) algorithm able to deal with large-scale instances, e.g. stochastic gradient descent (SGD).

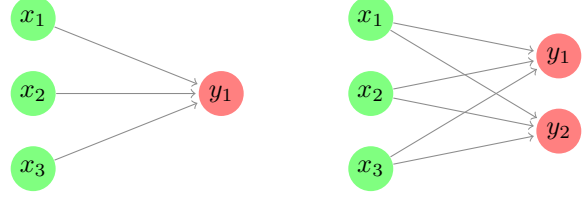

Figure 1: Standard (left) and multiple (right) logistic regressions.

## 3.2 PMLP-PaI

In this variant, we want to output an approximation $\bar{\mathcal{C}}(f, r)$ of the map that, given an instance and a desired performance in $[0, 1]$, returns the most appropriate configuration. Therefore, $X = F \times [0, 1]$ and $Y = C$ in the training set, i.e., $x = (f, r)$ is a pair (features, performance) and $y \in \mathbb{R}^s$ is a configuration. By the definition of $y$ in PMLP-PaI, the LR requires multiple output nodes instead of a single one (see Fig. 1 (right)), since $s > 1$ in general. This can simply be achieved by considering $s$ standard LRs sharing the same input node set.

**Proposition 3.1.** *The training problem of a multiple LR with $k$ output nodes consists of $k$ independent training problems for standard LRs, as in Eq. (7).*

*Proof.* A multiple LR on $k$ outputs is equivalent to $k$ standard LRs with training sets $T^1 = (X, Y^1), \ldots, T^k = (X, Y^k)$ where $Y^h = (y_h^1, \ldots, y_h^n)$ for all $h \leq k$ and $t \leq n$. Note that all these training sets share the same input vector set $X$. For each $h \leq k$ we define Bernoulli random variables $\mathcal{Y}_h$. Then $\mathsf{P}(\mathcal{Y}_h = 1 \,|\, \mathcal{X} = x)$ (for some $x \in \mathbb{R}^m$) is given by $\tau(x, w^h, b^h)$, where $w^h \in \mathbb{R}^m$ and $b^h \in \mathbb{R}$, for all $h \leq k$. The training problem aims at maximizing the log-likelihood functions $\ln L_{T^h}(w^h, b^h)$ Eq. (7) of each output node $h \leq k$, which yields the objective function $\max \sum_{h \leq k} \ln L_{T^h}(w^h, b^h)$. Now we note that the optimum of $\sum_h \ln L_{T^h}(w^h, b^h)$ is achieved by optimizing each term separately, since each term depends on separate decision variables. □

As anticipated, it is already rather hard to have the LR to produce a bona fide $y \in \{0, 1\}^s$, although this might be easily solved by rounding; what is much harder to obtain is that $y \in C$. Since $s > 1$ in general, the LR requires multiple output nodes instead of one f

## 4 The CSSP

The CSSP is the problem of computing a good configuration $c^*$ for an input instance $\bar{f}$ and the learnt PMLP map. Clearly, its formulation depends on the output of the learning phase, that is, either $\bar{p}(f, c)$ or $\bar{\mathcal{C}}(f, r)$. However, the solution of both the PaO and the PaI variant is guaranteed to be feasible w.r.t. all the dependence/compatibility constraints, i.e. $c^* \in \{c \,|\, c \in \{0, 1\}^s, Ac \leq d\}$.

### 4.1 CSSP-PaO

In this case, the most obvious version of the CSSP would be to just maximize the expected performance over the set of feasible configurations, consistently with the dynamics of the trained LR. This yields the following (nonconvex) Mixed-Integer NonLinear Program (MINLP):

$$\max \left\{ \frac{1}{1 + e^{-(w^*(\bar{f}, c) + b^*)}} \,\Big|\, Ac \leq d \,,\, c \in \{0, 1\}^s \right\}, \qquad (8)$$

for given $\bar{f}$. Note that Eq. (8) depends on the instance at hand through the input parameters $\bar{f}$, $w^* \in \mathbb{R}^{t+s}$ and $b^* \in \mathbb{R}$. As already remarked in Sect. 2.2, the objective function of (8) is log-concave, which means that

$$\max \left\{ \ln\left(\frac{1}{1 + e^{-(w^*(\bar{f}, c) + b^*)}}\right) \,\bigg|\, Ac \leq d \,,\, c \in \{0, 1\}^s \right\} \qquad (9)$$

is a MINLP yielding the same optima as (8).

We identified a different interpretation for the CSSP objective, namely that of maximizing the likelihood that any new instance would be matched with a solver configuration and a performance value "as closely as possible" to the associations between $(f, c)$ and $p(f, c)$ established during training. In other words, we maximize the likelihood given in (6) as a function of $c$ and $r$, $r$ being a specific performance value. In order to have the CSSP pick out a high performance, we add a term $+r$ to the objective:

$$\max_{c,r} r \ln\left(\frac{1}{1+e^{-(w^*(\bar{f},c)+b^*)}}\right) + (1-r)\left(1 - \ln\left(\frac{1}{1+e^{-(w^*(\bar{f},c)+b^*)}}\right)\right) + r \quad (10)$$

$$Ac \leq d \ , \ c \in \{0,1\}^s \ , \ r \in [0,1] \quad (11)$$

Note that, while the performance measure $r$ is not binary, it is in $[0,1]$ (where 1 corresponds to maximum (excellent) performance), which is compatible with LR.

Finally, we tested a third CSSP interpretation, where each alternative $r$ and $1-r$ is weighted by the corresponding conditional probability:

$$\max_{c,r} \left\{ r\left(\frac{1}{1+e^{-(w^*(\bar{f},c)+b^*)}}\right) + (1-r)\left(1 - \left(\frac{1}{1+e^{-(w^*(\bar{f},c)+b^*)}}\right)\right) \mid (11) \right\} \quad (12)$$

While (12) is non-convex, we were still able to (heuristically) solve it efficently enough using `bonmin`.

### 4.2 CSSP-PaI

We now consider the multiple LR setting which correlates a given instance feature/performance vector $(f,p)$ to a configuration $c$. Although $p$ being part of the input means we need not restrict it to $[0,1]$, we chose to replace it with a ranked and scaled version $r$ for better comparing with CSSP-PaO. The most direct interpretation of the CSSP in this case is the nonconvex MINLP

$$\max_{c,r} \left\{ r \mid c_j = \frac{1}{1+e^{-\langle (w^j)^*, (\bar{f},r)\rangle - b_j^*}} \ \forall j \leq s \ , \ (11) \right\} , \quad (13)$$

where $(w^j)^* \in \mathbb{R}^{t+1}$ is the weight vector of the $j$-th output. However, this interpretation does not satisfy the feasibility requirements on the $c_j$.

**Proposition 4.1.** *(13) is infeasible, even if the constraint $r \in [0,1]$ is relaxed.*

*Proof.* The constraint of the problem implies for all $j \leq s$

$$\frac{1}{1+e^{-\langle (w^j)^*, (\bar{f},r)\rangle - b_j^*}} \in \{0,1\}.$$

However, for any given $(w^j)^*$, $b_j^*$ and $\bar{f}$, there is no value of $r \in \mathbb{R}$ which makes the LHS either 0 or 1, hence the result. □ □

Because of Prop. 4.1, we consider the same interpretation of the CSSP yielding the best objective function for the PaO case, i.e., the MINLP

$$\max_{c,r} \left\{ \sum_{j \leq s} \left[ c_j \left(\frac{1}{1+e^{-\langle (w^j)^*, (\bar{f},r)\rangle - b_j^*}}\right) + (1-c_j)\left(1 - \left(\frac{1}{1+e^{-\langle (w^j)^*, (\bar{f},r)\rangle - b_j^*}}\right)\right) \right] \mid (11) \right\}$$

which, through simple rearrangements, can be reformulated as

$$\max_{c,r} \left\{ \sum_{j \leq s} \frac{\left(1 - e^{-\langle (w^j)^*, (\bar{f},r)\rangle - b_j^*}\right)c_j - 1}{1+e^{-\langle (w^j)^*, (\bar{f},r)\rangle - b_j^*}} \mid (11) \right\} \quad (14)$$

## 5 Computational experiments

We tested both the PaO and PaI variants of our approach in the following general set-up:

- we consider 41 mixed-integer linear programming instances of the Hydro Unit Commitment (HUC) problem [5];
- the MP solver of choice is `CPLEX` [17];
- the supervised ML methodology used in the PMLP is LR [7, 28];
- the CSSP is a MINLP which we heuristically solve—using the `bonmin` open-source solver [4]—to find good parameter values for CPLEX deployed on 41 instances of the HUC problem.

## 5.1 Technical specifications

All experiments were carried out on a single virtual core of a 1.4GHz Intel Core i7 of a MacBook 2017 with 16GB RAM running under macOS Mojave 10.14.6. Our implementations are based on Python 3.7 [27], AMPL 20200430 [12], and bonmin 1.8.6 [4]. We implemented LR as a Keras+TensorFlow [10, 1] neural network with sigmoid activation and a stochastic gradient descent solver minimizing a loss function given by binary cross-entropy (a simple reformulation of the log-likelihood function (7)). The ranking function turning the performance data into $\rho$ was supplied by `scipy.stats.rankdata` [19], and the scaling to $[0, 1]$ by `sklearn.preprocessing.minmax_scale` [24].

## 5.2 The algorithmic framework

In this section we give a detailed description of the general algorithmic framework we employ.

1. *Feature extraction.* A set of $t = 54$ features was extracted from each of the 41 problem instances, so $|F| = 41$ and $F \subseteq \mathbb{R}^{54}$.

2. *Selection of configuration parameters.* We considered a subset of 11 CPLEX parameters (`fpheur`, `rinsheur`, `dive`, `probe`, `heuristicfreq`, `startalgo-rithm`, `subalgorithm` from `mip.strategy`; `crossover` from `barrier`; and `mircuts`, `flowcovers`, `pathcut` from `mip.cuts`), each with a varying number of discrete settings (between 2 and 4), which we combined so as to obtain 9216 configurations. We transformed each of these settings into binary form, obtaining $s = 27$ binary parameters, so $|C| = 27$. These parameters were chosen because, in our experience, they were reasonably likely to have an impact on the problem we considered. Therefore, our dataset is composed of $41 \times 9216 = 377856$ points.

3. *Obtaining the performance data.* For each $(f, c) \in F \times C$ in the dataset, we ran three times, with different random seeds, CPLEX configured by $c$ over the instance described by $f$ for 60 seconds, recording as $p(f, c)$ the second best integrality gap attained (the closer to zero, the better); this allows to mitigate the effect of performance variability issues, by which MIP solvers such as CPLEX have been shown to be affected [21]. We then form the performance value list $\mathfrak{p} = (p(f, c) \mid (f, c) \in F \times C)$.

4. *Ranking and scaling.* We ranked $\mathfrak{p}$, scaled it to $\hat{\mathfrak{p}} \subset [0, 1]$, and let $\rho = \mathbf{1} - \hat{\mathfrak{p}}$ in order for the value 1 to mean "best performance".

5. *Separating in-sample and out-of-sample sets.* We randomly choose 11 out of the 41 instances as "out-of-sample", put them in a set $F''$, and let $F' = F \smallsetminus F''$ be the "in-sample" set.

6. *Construction of the training sets.* For PMLP-PaO we let $X = F' \times C$ and $Y = \rho$, while for PMLP-PaI we let $X = F' \times \rho$ and $Y = C$.

7. We use the `sklearn.cluster.KMeans` k-means algorithm implementation to cluster the dataset into 5 clusters. We form a training set with 75% of the vectors from each cluster, a validation set with 20%, and a test set with the remaining 5%. By using clustering, we want to ensure that, even after the sampling, the actual distribution of the instances is preserved in all sets.

8. We implement a LR using a `keras.layers.Dense` complete bipartite pair of input/output layers (for PMLP-PaO with $f + s$ input nodes and 1 output node, for PMLP-PaI with $f + 1$ input nodes and $s$ output nodes), with a sigmoid activation function in the output nodes. We train the LR using the corresponding training, validation and test sets, using the keras SGD optimizer optimizing the binary cross-entropy loss function (which corresponds to minimizing the negative of Eq. (7)). Then, for further use in the different CSSP (as described in Sect. 4), we save:
    - $(w^*, b^*)$, with $w \in \mathbb{R}^{t+s}$ and $b \in \mathbb{R}$, for PMLP-PaO;
    - $((w^j)^*, b_j^*)$, $\forall j \leq s$, with $w^j \in \mathbb{R}^{t+1}$ and $b_j \in \mathbb{R}$, for PMLP-PaI.

9. For each out-of-sample instance feature vector $g \in F''$ we perform the following actions:
    (a) we establish a link from Python to AMPL via `amplpy`;
    (b) we solve the CSSP corresponding to the feature vector $g$ with `bonmin`;
    (c) we retrieve the optimal configuration $c^*$;
    (d) we retrieve the stored performances $p(g, c^*)$ and $p(g, d)$, where $d$ is the default CPLEX configuration;
    (e) if $p(g, c^*) > p(g, d)$ we count an improvement;
    (f) if $p(g, c^*) \geq p(g, d) - 0.001$ we count a non-worsening;

(g) we record the performance difference $|p(g, c^*) - p(g, d)|$.

10. We count the number of improvements im and non-worsenings nw over the number of successful `bonmin` runs on the CSSP instances; sometimes `bonmin` fails on account of the underlying NLP solver, which is why some lines of Table 1 consider a total of less than 11 instances.

11. We repeat this process 10 times from Step 5, and report cumulative statistics of improvements, non-worsenings, performance differences, and CPU time.

## 5.3 Results

We first conducted experiments on the simple PaO interpretation (9) of the CSSP. However, this gave very poor results in practice. Eq. (11), instead, gave better computational results than those obtained optimizing (9), although each CSSP instance took considerably more time to solve w.r.t. (9) and (12). The PaO formulation (12) is the one which gave the best results and is therefore the only one considered in Table 1. As for the PaI variant, Table 1 shows the results of formulation (14). In the table below, we report improvements im, non-worsenings nw, performance differences pd, and CPU times. We also report cumulative statistics (sum, mean, standard deviation) for the 10 runs of the algorithmic framework in Sect. 5.2 for the PaO and PaI variants. We remark that the "by-run" comparison is only meant for presentation, as the out-of-sample instances involved in each run of PaO and PaI differ.

|       | im       |          | nw       |          | pd       |          | CPU      |          |
|-------|----------|----------|----------|----------|----------|----------|----------|----------|
| run   | PaO      | PaI      | PaO      | PaI      | PaO      | PaI      | PaO      | PaI      |
| 1     | 0/08     | **5/09** | 0/08     | **7/09** | 0.63     | **0.30** | 41.44    | **30.55** |
| 2     | 0/11     | **4/11** | 0/11     | **6/11** | **0.42** | 0.47     | 41.62    | **29.28** |
| 3     | 4/09     | 4/10     | 5/09     | 8/10     | **0.08** | 0.14     | 43.06    | **33.37** |
| 4     | 0/09     | **5/10** | 0/09     | **8/10** | 0.43     | **0.12** | 42.65    | **35.28** |
| 5     | **3/10** | 1/10     | **7/10** | 2/10     | **0.08** | 0.70     | 43.30    | **31.69** |
| 6     | **8/09** | 3/11     | **8/09** | 9/11     | 0.20     | **0.18** | 43.69    | **28.98** |
| 7     | **5/10** | 1/11     | **8/10** | 4/11     | **0.05** | 0.52     | 45.54    | **30.05** |
| 8     | **7/08** | 3/09     | 7/08     | 8/09     | 0.21     | **0.02** | 45.49    | **31.28** |
| 9     | 0/09     | 0/10     | 0/09     | 0/10     | **0.40** | 0.88     | 43.83    | **31.88** |
| 10    | **8/10** | 5/08     | 8/10     | 7/08     | 0.21     | **0.10** | 43.50    | **33.88** |
| sum   | **35/93** | 34/99   | 43/93    | 59/99    | **2.69** | 3.40     | 434.12   | 316.24   |
| mean  | **0.38** | 0.32     | 0.46     | **0.60** | **0.26** | 0.34     | 43.41    | **31.62** |
| stdev | 0.36     | **0.20** | 0.39     | **0.30** | **0.18** | 0.27     | **1.30** | 1.94     |

Table 1: Computational results. Best results are marked in boldface.

The results show that the PaO and PaI variants are comparable. PaO improves more, but also worsens more. PaI improves slightly less, but it has three considerable advantages w.r.t. PaO: (i) it does not worsen results more than 60% of the times, which means it can be recommended for usage w.r.t. the default CPLEX configuration; (ii) it is more reliable in terms of standard deviation of improvements and non-worsening; (iii) it is faster.

## 5.4 Conclusions

We presented a general two-phase framework for learning good mathematical programming solver configurations, subject to logical constraints, using a performance function estimated from data. We proposed two significantly different variants of the methodology, both using Logistic Regression as the Machine Learning technique of choice, but using different configurations for the inputs and outputs of the LR. Tested on a problem arising in scheduling of hydro-electric generators, both variants showed promise, although the PaI one appeared to be preferable for several reasons. We remark that these encouraging results were obtained with a relatively small number of instances.

In future works, we are going to investigate our general framework using different Machine Learning methodologies, such as (deep) Neural Networks and Support Vector Machine/Regression. Each ML technique requires a different definition of both the PMLP and the CSSP (for each of the two PaO and PaI variants); hence, the exploration of the vast landscape of possible versions will offer a vast choice of the trade-offs between computational cost and effectiveness of the obtained configuration, hopefully finally leading to versions that may become actually useful for day-to-day use of MP solvers.